\documentclass[11pt]{amsart}
\usepackage{ytableau}
\usepackage{tikz}
\usepackage{mathtools}
\usepackage{enumitem}
\usepackage{fullpage}

 \newtheorem{theorem}{Theorem}[section]
 \newtheorem{corollary}[theorem]{Corollary}
 \newtheorem{lemma}[theorem]{Lemma}
 \newtheorem{proposition}[theorem]{Proposition}
 \theoremstyle{definition}
 \newtheorem{definition}[theorem]{Definition}
 \newtheorem{remark}[theorem]{Remark}

  \newtheorem{example}[theorem]{Example}
  
  \newtheorem{qu}[theorem]{Question}
 \numberwithin{equation}{section}
\DeclareMathOperator{\Ima}{Im}
\DeclareMathOperator{\Hom}{Hom}
\DeclareMathOperator{\coker}{Coker}
\DeclareMathOperator{\sgn}{sgn}

\begin{document}

\title{On the action of the symmetric group on the free LAnKe: a question of Friedmann, Hanlon, Stanley and Wachs}
%----------Author 1

\author[Maliakas]{Mihalis Maliakas}
\address{%
	Department of Mathematics\\
	University of Athens\\
	Greece}
\email{mmaliak@math.uoa.gr}

\author[Stergiopoulou]{Dimitra-Dionysia Stergiopoulou}
\address{School of Applied Mathematical and Physical Sciences, National Technical University of Athens, Greece}
	\curraddr{Department of Mathematics\\
		University of Athens\\
		Greece}
\email{dstergiop@math.uoa.gr}



\date{March 15, 2025}

%%% ----------------------------------------------------------------------

\begin{abstract}A LAnKe (also known as a Lie algebra of the $n$th kind, or a Filippov algebra) is a vector space equipped with a skew-symmetric $n$-linear form that satisfies the generalized Jacobi identity. The symmetric group $\mathfrak{S}_m$ acts on the multilinear part of the free LAnKe on $m=(n-1)k+1$ generators, where $k$ is the number of brackets, by permutation of the generators. The corresponding representation was studied by Friedmann, Hanlon, Stanley and Wachs, who asked whether for $n \ge k$, its irreducible decomposition contains no  summand whose Young diagram has at most $k-1$ columns. The answer is affirmative if $k \le 3$. In this paper, we show that the answer is affirmative for all $k$. A proof has been given recently by Friedmann, Hanlon and Wachs. The two proofs are completely different.
\end{abstract}
\maketitle
\tableofcontents

	\section{Introduction}Since the mid 1980's, various  $n$-ary generalizations of Lie algebras have been studied. These include the algebras introduced by Filippov \cite{Fi},  the Lie $n$-algebras introduced by Hanlon and Wachs \cite{HW} and the Liebniz $n$-algebras introduced by Casas, Loday and Pirashvili.   In the cases of the algebras introduced by Filippov and by Hanlon and Wachs, the representation of the symmetric group acting on the multilinear component of the free $n$-ary Lie algebra generalizes the Lie representation which has been extensively studied, for example see \cite{Ba, Re, Wa, Su}. We refer to the Introduction of the paper by Friedmann, Hanlon, Stanley and Wachs \cite{FHSW} for a relevant discussion including historical background and ties with various areas of mathematics and physics.
	
	We work over a field $\mathbb{K}$ of characteristic zero. Let $\mathfrak{S}_n$ be the symmetric group on $n$ symbols. A \textit{Lie algebra of the $n$th kind} (\textit{LAnKe}, or \textit{Filippov algebra}) is a $\mathbb{K}$-vector space $\mathcal{L}$ equipped with an $n$-linear bracket \[ [-,-, \dots, -] :  \mathcal{L}^n \to \mathcal{L}\] such that for all $x_1,\dots, x_n, y_1, \dots, y_{n-1} \in \mathcal{L}$,
	\begin{enumerate}
		\item $[x_1, \dots, x_n]=\sgn(\sigma)[x_{\sigma(1)}, \dots, x_{\sigma(n)}$] for every $\sigma \in \mathfrak{S}_n$, and
		\item the following generalized Jacobi identity holds\begin{align}\label{GJI}
			[[x_1,\dots,x_n],y_1,\dots, y_{n-1}]
			=\sum_{i=1}^{n}[x_1,\dots,x_{i-1},[x_i,y_1,\dots, y_{n-1}],x_{i+1},\dots, x_n].
	\end{align}	\end{enumerate}
	
	These algebras were introduced by Filippov \cite{Fi} and have attracted the attention of physicists, for example see \cite{AI, Fr}. More recently, Friedmann, Hanlon, Stanley and Wachs \cite{FHSW} initiated the study of the action of the symmetric group on the multilinear component of the free LAnKe.
	
	Let us consider the free LAnKe on $[m]:=\{1, \dots, m\}$. The free LAnKe for $n=2$ is the free Lie algebra. Following \cite{FHSW}, the \textit{multilinear component} $\text{Lie}_n(m)$ of the free LAnKe on $[m]$ is spanned by the bracketed words on $[m]$ in which each $i$ appears exactly once. It follows that each such bracketed word has the same number of brackets, say $k$, and $m=(n-1)k+1$. Consider the action of the symmetric group $\mathfrak{S}_{m}$ on $\text{Lie}_n(m)$ given by replacing $i$ by $\sigma(i)$ in each bracketed word. We denote the corresponding representation of $\mathfrak{S}_m$ by $\rho_{n,k}$. Since the characteristic of $\mathbb{K}$ is equal to 0, the irreducible representations of $\mathfrak{S}_{m}$ are indexed by the partitions $\lambda$ of $m$. We denote by $S^{\lambda}$ the Specht module corresponding to $\lambda$. This is an irreducible $\mathfrak{S}_{m}$-module and as $\lambda$ runs over the partitions of $m$, the Specht modules form a complete set of irreducible $\mathfrak{S}_{m}$-modules.
	\subsection{Motivation} The irreducible decomposition of $\rho_{n,k}$ remains an open problem. Let us recall an important result of Friedmann, Hanlon, Stanley and Wachs from \cite{FHSW} and \cite{FHW2} that is the main motivation of the present paper. 
	
	For $n \ge 2$, let $\beta_{n,k}$ be the $\mathfrak{S}_m$-module whose decomposition into irreducibles is obtained by adding a row of length $k$ to the top of each Young diagram appearing in the decomposition into irreducibles of $\rho_{n-1, k}$. For $n=1$ we regard $\rho_{1,k}=\beta_{1,k}=S^{(1)}$.
	\begin{theorem}{\emph{\cite[Theorem 1.5]{FHW2}}} \label{dec} Let $n, k \ge 1$. Then, as $\mathfrak{S}_m$-modules, \begin{equation}\label{FHWdec}
			\rho_{n,k} \simeq \beta_{n,k} \oplus \gamma_{n,k}\end{equation}
		for some $\mathfrak{S}_m$-module $\gamma_{n,k}$ all of whose irreducibles have Young diagrams with at most $k-1$ columns. 
	\end{theorem}	
	The following question is one of the central open problems posed by Friedmann, Hanlon, Stanley and Wachs in \cite{FHSW}.
	\begin{qu}{\cite[Question 4.1]{FHSW}} \label{qu} Does $n \ge k$ imply $\gamma_{n,k}=0$?
	\end{qu}
	It may happen that $\gamma_{n,k} \neq 0$ if $n<k$. For example,  $\gamma_{3,4} \neq 0$ \cite[proof of Theorem 7.1]{FHW2}. 
	
	It was shown in \cite[Theorem 1.3]{FHSW} that $\rho_{n,2}$ is an irreducible module and thus Question \ref{qu} has an affirmative answer when $k=2$. 
	
	The decomposition of $\rho_{n,3}$ into irreducibles was announced in \cite{FHSW}. For two completely different proofs see \cite{MS7} and \cite[Theorem 1.3]{FHW2}. This result implies that Question \ref{qu} has an affirmative answer for $k=3$ as well.
	
	\subsection{Results} In this paper, we show that the answer to Question \ref{qu} is affirmative for all $n \ge k$, see Theorem \ref{main}(2). From this and Theorem \ref{dec} we see that, for $k$ fixed, the irreducible decompositions of the representations $\rho_{k-1,k}$, $\rho_{k,k}$, $\rho_{k+1,k}$, $\dots$ enjoy a pleasing recursion.
	
	The main idea of the proof is to show that for all $n \ge 2$ and $k \ge 3$ the $\mathfrak{S}_m$-module $\mathrm{Lie}_n(m)$ is a homomorphic image of the skew Specht module $S^{\beta(n,k)}$, where $\beta(n,k)$ is the skew partition
	\begin{equation}\label{beta}\beta(n,k):= (k^{n-1}, k-1, k-2, \dots, 1)/(k-1, k-2, \dots, 2),\end{equation} where $k^{n-1}$ means that the part $k$ appears $n-1$ times.
	\begin{center}
		\ytableausetup{boxsize=1.2em}
		\ydiagram{3+1,2+2,4,3,2,1} \ \ \ \ \ \ytableausetup{boxsize=1.2em}
		\ydiagram{3+1,2+2,3,2,1} \\ \; \\
		$\beta(4,4)$ \ \ \ \ \ \ \ \ \ \ \ $\beta(3,4)$
	\end{center}
	
	Now if $n \ge k$, the Littlewood-Richardson rule implies that the multiplicity in $S^{\beta(n,k)}$ of any Specht module whose Young diagram has at most $k-1$ columns is equal to zero and hence we have $\gamma_{n,k} =0$.
	
	The decomposition of $\rho_{n,k}$ into irreducibles is known for $k \le 3$, see \cite{FHSW, FHW2, MS7}. As a corollary of our main result and \cite[Theorem 7.1]{FHW2}, we obtain the decomposition of $\rho_{n,4}$, see Corollary \ref{cordec}. We remarked previously that, in general, the decomposition of $\rho_{n,k}$ remains an open problem. Our main result though, yields a new upper bound on the multiplicities of the irreducibles valid for all ${n,k}$, see Corollary \ref{cordecsp}.
	\begin{remark} The first version of the present paper (arXiv:2410.06979) was posted on October 10, 2024. A proof of an affirmative answer to Question \ref{qu}  was added to the third version (arXiv:2402.19174v3) of \cite{FHW2} posted on October 13, 2024. This proof appears to be completely different from ours. 
	\end{remark}
	This paper is organized as follows. We approach the problem from the point of view of representations of the general linear group. Relevant preliminaries are gathered in Section 2. In Section 3 we study specific relations of the multilinear component $\mathrm{Lie}_n{(m)}$ that are central for our purposes. In Section 3 and Section 4, the $GL_N$ versions of these relations are shown to coincide with the defining relations of certain skew Weyl modules. Then, we prove our main result on Weyl modules in Section 5 (Theorem \ref{mainweyl}). In Section 6, we use the duality functor $\Omega$ and the Schur functor to derive an affirmative  answer of Question \ref{qu} from Theorem \ref{mainweyl}.
	\section{Preliminaries}
	The purpose of this section is to gather material that will be used in the sequel. This concerns mainly the divided power algebra and the exterior algebra, skew Weyl modules and combinatorics of tableaux, the duality functor $\Omega$ and the Schur functor $f$.\subsection{Divided power algebra and exterior algebra}
	Let $G=GL_N(\mathbb{K})$ be the general linear group of $N \times N$ matrices with entries in $\mathbb{K}$. Let $V=\mathbb{K}^N$ be the natural $G$-module consisting of column vectors.  By $D=\sum_{i\geq 0}D_i$ we denote the divided power algebra of $V$ \cite[Section 1.1]{W}. By definition we have $D_i= (S_i(V^*))^*$, where $S_i(V^*)$ denotes the degree $i$ symmetric power of the dual $V^*$ of the vector space $V$. Since the characteristic of $\mathbb{K}$ is equal to 0, the divided power algebra of  $V$ and the symmetric algebra of $V$ are naturally isomorphic. However, we will work with the divided power algebra since the computations of Section 4 below seem less involved.  
	
	If $v \in V$ and $i$ is a nonnegative integer, we have the $i$th divided power $v^{(r)} \in D_i$. We recall that if $i,j$ are nonnegative integers, then the product $v^{(i)}v^{(j)}$ of $v^{(i)}$ and $v^{(j)}$ is given by \[v^{(i)}v^{(j)}=\tbinom{i+j}{j}v^{(i+j)},\] where $\tbinom{i+j}{j}$ is the indicated binomial coefficient. 
	
	If $\{e_1, \dots, e_N\}$ is a basis of $V$, then a basis of $D_i$ is the set \[\{e_1^{(\alpha_1)}\cdots e_N^{(\alpha_N)}: \alpha_1+\cdots \alpha_N = i.\}\]
	We recall that $D$ has a graded Hopf algebra structure. Let \[\Delta : D \to D \otimes D\] be the comultiplication map of $D$. Explicitly, for a homogeneous element $x= v_1^{(\alpha_1)}\cdots v_t^{(\alpha_t)} \in D_a$, where $v_i \in V$, we have \[ \Delta(x)=\sum_{0\le \beta_i \le \alpha_i} v_1^{(\beta_1)}\cdots v_t^{(\beta_t)} \otimes v_1^{(\alpha_1 - \beta_1)}\cdots v_t^{(\alpha_t - \beta_t)}.\]
	For $0 \le b \le a$ we may restrict the above sum to those $\beta_i$ such that $\beta_1 + \cdots +  \beta_t= b$. This yields the following component of the comultiplication map \begin{align*} D_a &\to D_b \otimes D_{a-b}, \\  x &\mapsto \sum_{\substack{0\le \beta_i \le \alpha_i \\ \beta_1+\cdots+\beta_t=b}} v_1^{(\beta_1)}\cdots v_t^{(\beta_t)} \otimes v_1^{(\alpha_1 - \beta_1)}\cdots v_t^{(\alpha_t - \beta_t)}, \end{align*}
	which we will again denote simply by $\Delta:  D_a \to D_b \otimes D_{a-b}$ in order to avoid cumbersome notation.

	By $\Lambda = \bigoplus_{i \ge 0} \Lambda^i$ we denote the exterior algebra of $V$. We recall that $\Lambda$ has a graded Hopf algebra structure. If $u,v \in \Lambda$, we denote their product in  $\Lambda$ by $uv$. If $\{e_1, \dots, e_N\}$ is a basis of the vector space $V$, then a basis of the vector space $\Lambda^i$ is the set \[\{e_{\alpha_1} e_{\alpha_{2}}\cdots  e_{\alpha_i}: 1 \le \alpha_1 < \cdots < \alpha_i \le N\}.\] We denote by \[\Delta : \Lambda \to \Lambda \otimes \Lambda\] the comultiplication map of $\Lambda$. Explicitly, for a homogeneous element $x= v_{1} v_{2}\cdots  v_{a} \in \Lambda^a$, where $v_i \in V$, we have \[ \Delta(x)=\sum_{0 \le s \le a}\sum_{\sigma} v_{\sigma(1)} \dots v_{\sigma(s)} \otimes v_{\sigma(s+1)} \dots v_{\sigma(a)},\] where the second sum is over all permutations $\sigma$ of $\{1, \dots, a\}$ such that $\sigma(1) < \dots <\sigma(s)$ and $\sigma(s+1) < \dots < \sigma(a)$.
	
	For $0 \le b \le a$ we have the following component of the comultiplication map \begin{align*} \Lambda^a &\to \Lambda^b \otimes \Lambda^{a-b}, \\  x &\mapsto \sum_{\sigma} v_{\sigma(1)} \dots v_{\sigma(b)} \otimes v_{\sigma(s+1)} \dots v_{\sigma(a)}, \end{align*} where the sum is over all permutations $\sigma$ of $\{1, \dots, a\}$ such that $\sigma(1) < \dots <\sigma(s)$ and $\sigma(s+1) < \dots < \sigma(a)$. We will  denote this map simply by $\Delta:  \Lambda^a \to \Lambda^b \otimes \Lambda^{a-b}$.
	
	We have used the same symbol $\Delta$ for the comultiplication maps in the algebras $D$ and $\Lambda$. In the sequel it will be clear which algebra is considered each time. If there is a need of distinction, we will write $\Delta_D$ and $\Delta_\Lambda$.

	\subsection{Skew partitions, skew Weyl modules and skew Schur modules}For a positive integer $r$, a \textit{partition} of $r$ is a sequence  $\lambda=(\lambda_1, \dots, \lambda_q)$ of nonnegative integers such that $\lambda_1 \ge \dots \ge \lambda_q$ and $\lambda_1  +\dots + \lambda_q=r$. If $\lambda$ is a partition of $r$ we may write $r=|\lambda|$. We will identify two partitions if they differ by a string of zeros at the end. The conjugate of a partition $\lambda$ is denoted by $\lambda '$.
	
	A \textit{skew partition} $\lambda/\mu$ is a pair of partitions $\lambda =(\lambda_1, \dots, \lambda_q)$, $\mu =(\mu_1, \dots, \mu_q)$ such that $\lambda_ i \ge \mu_i$ for all $i$.  The Young diagram of a skew partition can be obtained from the Young diagram of $\lambda$ by omitting the cells that correspond to $\mu$. We use the English convention for diagrams of partitions. For example, the Young diagram of the skew partition $(6,4,3)/(1,1)$ is
	\begin{center}
		\ytableausetup{smalltableaux}
		\ydiagram{1+5,1+3,3}
	\end{center}
	Often we will identify a skew partition with its Young diagram.
	
	For a skew partition $\lambda/\mu$, we denote by $K_{\lambda/\mu}$ the corresponding skew Weyl module for $G$ \cite[Definition II.1.4]{ABW}. (In loc. cit. the term `CoSchur functor' is used for skew Weyl modules). The formal character of $K_{\lambda/\mu}$ is the skew Schur polynomial $s_{\lambda/\mu}$. When $\mu=(0)$ we will write $K_\lambda$ in place of $K_{\lambda / \mu}$. For example, when $\lambda=(r)$ consists of one part and $\mu=(0)$, then $K_{\lambda} = D_r$ and the formal character is the complete symmetric polynomial of degree $r$.
	
	For a skew partition $\lambda/\mu$, we denote by $L_{\lambda/\mu}$ the corresponding skew Schur module for $G$ \cite[Definition II.1.3]{ABW}. The formal character of $L_{\lambda/\mu}$ is the skew Schur polynomial $s_{\lambda'/\mu'}$. When $\mu=(0)$ we will write $L_\lambda$ in place of $L_{\lambda / \mu}$. For example, when $\lambda=(r)$ consists of one part and $\mu=(0)$, then $L_{\lambda} = \Lambda^r$ and the formal character is the elementary symmetric polynomial of degree $r$.
	
	For a sequence $\alpha=(\alpha_1, \cdots, \alpha_q)$ of nonnegative integers we let  \begin{align*}&D(\alpha):=D_{\alpha_1} \otimes \dots \otimes D_{\alpha_q},\\&\Lambda(\alpha):=\Lambda^{\alpha_1} \otimes \dots \otimes \Lambda^{\alpha_q}\end{align*} be the tensor product of the indicated of divided powers, respectively exterior powers, over $\mathbb{K}$. 
 	
	If $\lambda$ is a partition, then there is a unique (up to scalar) projection of $G$-modules $D(\lambda) \to K_\lambda$ which we denote by $\pi_\lambda$,\[ \pi_\lambda : D(\lambda) \to K_\lambda.\]
	
	Suppose $N \ge r$. Since the characteristic of $\mathbb{K}$ is zero, the $G$-module $K_\lambda$ is irreducible for every partition $\lambda$ of $r$ and, moreover, for distinct partitions $\lambda$ and $\mu$ of $r$, the $G$-modules $K_\lambda$ and $K_\mu$ are non isomorphic. Every polynomial representation of $G$ is a direct sum of various $K_\lambda$. See \cite[Section 2.2]{W}.
	
	The skew Weyl module $K_{\lambda / \mu}$ is not in general irreducible and its irreducible decomposition is given by  \[K_{\lambda / \mu} = \sum_{\nu} c^{\lambda}_{\mu, \nu}K_\nu,\]
	where the sum ranges over all partition $\nu$ such that $|\nu| = |\lambda| - |\mu|$. Here, $c^{\lambda}_{\mu, \nu}$ denotes the multiplicity of $K_\lambda$ as a summand of the tensor product $ K_\mu \otimes K_\nu$. The coefficients have a beautiful combinatorial description known as the  Littlewood-Richardson rule \cite[Section 5.2]{F}. 
	
	\subsection{Tableaux and maps}\label{2.2}
	
	Let us fix the order $e_1<e_2< \cdots <e_N$ on the natural basis $\{e_1,...,e_N\}$ of $V$. In the sequel we will denote  $e_i$ by its subscript $i$. If $\lambda/\mu$ is a skew partition, a \textit{tableau} of shape $\lambda/\mu$ is a filling of the Young  diagram of $\lambda/\mu$ with entries from $\{1,...,N\}$. A tableau is called \textit{row semistandard} if the entries are weakly increasing across the rows from left to right.  A row semistandard tableau is called \textit{semistandard} if the entries are strictly increasing down each column. We denote the set of row semistandard (respectively, semistandard) tableaux of shape $\lambda/\mu$ by $\mathrm{RSST}(\lambda/\mu)$ (respectively, $\mathrm{SST}(\lambda/\mu)$).
	
	We will use `exponential' notation for row semistandard tableaux corresponding to partitions. For example, we write \[\begin{matrix*}[l]
		1^{(3)} 2^{(2)} \\
		2^{(2)} 4^{(2)} \end{matrix*} \]
	for the tableau \begin{center}
		\begin{ytableau}
			1&1&1&2&2\\2&2&4&4
		\end{ytableau}
	\end{center}
	In the sequel we will need to define various $G$-maps of the form $D(a,b) \to D(c,d)$, where $a,b,c,d$ are nonnegative integers such that $a+b=c+d$. One way of describing such maps is with the use of two-rowed row standard tableaux. 
	\begin{definition}\label{phiS}Suppose 
		$ S:=\begin{matrix*}[l]
			1^{(a_{1})} 2^{(b_{1})} \\
			1^{(a_{2})} 2^{(b_{2})} \end{matrix*}.$ We define a map of $G$-modules \[\phi_S : D(a_1+a_2,b_1+b_2)  \to D{(a_1+b_1, a_2+b_2)}\] as the composition
		\begin{align}\label{phis} D(a_1+a_2,b_1+b_2) \xrightarrow{\Delta_1\otimes \Delta_2}& D(a_{1},a_{2}) \otimes D(b_1,b_2) \simeq D(a_{1},b_{1}) \otimes D(a_{2},b_2)  \\\nonumber\xrightarrow{\eta_1 \otimes \eta_2}& D(a_1+b_1, a_2+b_2),
		\end{align}
		where $\Delta_1 : D(a_1+a_2) \to D(a_1, a_2)$ and $\Delta_2 : D(b_1+b_2) \to D(b_1, b_2)$ are the indicated components of twofold comultiplication of the Hopf algebra $D$,  the isomorphism permutes tensor factors, and $\eta_i $: $D(a_{i}, b_{i}) \to D(a_i+b_i)$ are the indicated components of multiplication in the algebra $D$.\end{definition}
	
	\begin{example} For $S:=\begin{matrix*}[l]
			1^{(2)} 2^{(3)} \\
			12^{(3)}  \end{matrix*}$ we have the map \[\phi_S: D(3,6) \to D(5,4).\] Let $x=12^{(2)}\otimes 12^{(5)} \in D(3,6)$.
		With the notation of the previous definition  we have \begin{align*}
			&\Delta_1(12^{(2)})=12 \otimes 2 + 2^{(2)}\otimes 1, \\
			&\Delta_2(12^{(5)})=12^{(2)} \otimes 2^{(3)}+2^{(3)} \otimes 12^{(2)}.\end{align*} Hence the image of $x$ is 
		\begin{align*}
			\phi_S(x)&= 1122^{(2)}\otimes 22^{(3)}+122^{(3)}\otimes 122^{(2)} +12^{(2)}2^{(2)}\otimes 12^{(3)} +2^{(2)}2^{(3)}\otimes 112^{(2)}
			\\&=\tbinom{1+1}{1} \tbinom{1+2}{1}\tbinom{1+3}{1} 1^{(2)}2^{(3)} \otimes  2^{(4)} \\&+ \Big( \tbinom{1+3}{1}\tbinom{1+2}{1}+\tbinom{2+2}{2} \Big) 12^{(4)}\otimes 12^{(3)} +\tbinom{2+3}{2}\tbinom{1+1}{1} 2^{(5)}\otimes 1^{(2)}2^{(2)}.
		\end{align*}
		The binomial coefficients come from the multiplication in the algebra $D$.\end{example}
	
	\subsection{A presentation of skew Weyl modules}We will use the well known  presentation of skew Weyl modules described below. \begin{definition}\label{thetat}
		For nonnegative integers $a,b,t$ define the map \begin{equation*}\theta_t : D(a,b) \xrightarrow{\Delta \otimes 1} D(a-t, t,b) \xrightarrow{1 \otimes \eta}D(a-t,b+t),\end{equation*} where \[\Delta: D(a) \to D(a-t)\otimes D(t)\] is the indicated component of twofold comultiplication of the Hopf algebra $D$ and \[\eta: D(t,b) \to D(b+t)\] is the indicated component of multiplication.\end{definition}
	
	The relations of $K_{\lambda/\mu}$ described in the theorem that follows, are parametrized by pairs of consecutive rows of $\lambda/\mu$. Let us consider first the case of two rows. 
	
	Suppose $\lambda / \mu$ is a skew partition, $\lambda=(\lambda_1,\lambda_2)$, $\mu=(\mu_1,\mu_2)$. Then, as $G$-modules, $K_{\lambda/\mu}$ is isomorphic to the cokernel of the map  \begin{equation}\label{skew2} \theta_{\lambda/\mu} : D(\lambda_1-\mu_2+1, \lambda_2-\mu_1-1 ) \to D(\lambda / \mu),
	\end{equation}
	where $\theta_{\lambda/\mu} := \theta_{\mu_1-\mu_2+1}$. 
	
	In general we have the following presentation of $K_{\lambda/\mu}$.
	\begin{theorem}\label{skewmany}Let $\lambda/\mu$ be a skew partition, where $\lambda=(\lambda_1, \dots,\lambda_q)$, $\mu =(\mu_1, \dots, \mu_q)$. Let $\lambda ^{s}:=(\lambda_{s}, \lambda_{s+1})$ and $\mu ^{s}:=(\mu_{s}, \mu_{s+1})$, $s=1, \dots, q-1$. Then, as $G$-modules, $K_{\lambda/\mu}$ is isomorphic to the cokernel of the map \begin{equation}
			\theta_{\lambda/\mu}: \bigoplus_{s=1}^{q-1} D(\lambda_1-\mu_1, \dots,\lambda_s - \mu_{s+1}+1, \lambda_{s+1} - \mu_{s}-1, \dots, \lambda_q-\mu_q ) \to D(\lambda / \mu),
		\end{equation}where \[\theta_{\lambda/\mu}:=\sum_{s=1}^{q-1}1\otimes \cdots \otimes 1 \otimes \theta_{\lambda^s /\mu^s}\otimes 1 \otimes \cdots \otimes 1.\]
	\end{theorem} This follows from \cite[Theorem II.3.16]{ABW} and the fact that the characteristic of $\mathbb{K}$ is zero.
	\begin{example}For the skew partitions \[ \alpha = 
		\ydiagram{1+5,1+3}, \ \ \beta = 
		\ydiagram{1+3,3}, \ \ \gamma = 
		\ydiagram{1+5,1+3,3} \]
		we have the following descriptions of the corresponding Weyl modules	
		\begin{align*} &K_{\alpha}=\coker(\theta_1 : D(6,2) \to D(5,3)), \\
			&K_{\beta}=\coker(\theta_2 : D(5,1) \to D(3,3)), \\
			&K_{\gamma}=\coker(\theta_1 \otimes 1 \  +  \ 1 \otimes \theta_2: D(6,2,3) \oplus D(5,5,1) \to D(5,3,3)). 
		\end{align*}
		
	\end{example}
	\subsection{Straightening row semistandard tableaux} If $\nu = (\nu_1, \nu_2)$ is a partition with at most two parts and $S$ is a row semistandard tableau of shape $\nu$,
	\[ S=\begin{matrix*}[l]
		1^{(a_1)}2^{(a_2)}  \cdots   N^{(a_N)} \\
		1^{(b_1)}2^{(b_2)}  \cdots  N^{(b_N)}
	\end{matrix*},\]
	where $a_{i}, b_i$ are nonnegative integers, let $e^S \in D(\mu)$ be the element \[e^S:=1^{(a_{1})} \cdots N^{(a_{N})}  \otimes 1^{(b_{1})} \cdots N^{(b_{N})} \] obtained by `reading the rows' of $S$ from left to right and top to bottom.
	
	We recall a classical result (which is stated here only for partitions with at most two parts), see \cite[(2.1.15) Proposition]{W}.
	\begin{theorem}\label{sbthm} $\nu = (\nu_1, \nu_2)$ is a partition. Then there is a bijection between $\mathrm{SST}(\nu)$, and a basis of the $\mathbb{K}$-vector space $K_{\nu}$ given by 
			$ S \mapsto \pi_{\nu}(e^S).$\end{theorem}
		We refer to the elements of this basis of $K_\nu$ as \textit{semistandard basis elements}.
	
	In the sequel we will need to express elements of Weyl modules as explicit linear combinations of semistandard basis elements.   The following lemma concerns violations of semistandardness in the first column.
	
	\begin{lemma}[{\cite[Lemma 4.2]{MS3}}]\label{lemglas}Let $\nu=(\nu_1,\nu_2)$ be a partition of length two and let \[S=\begin{matrix*}[l]
			1^{(a_1)}2^{(a_2)}  \cdots   N^{(a_N)} \\
			1^{(b_1)}2^{(b_2)}  \cdots  N^{(b_N)}
		\end{matrix*} \in \mathrm{RSST} (\nu).\]
		Then we have the following identities in $K_{\nu}$.
		\begin{enumerate}
			\item If $a_1+b_1>\nu_1$, then $\pi_\nu(e^S)=0$.
			\item If $a_1+b_1 \le \nu_1$, then 
			\begin{equation}\label{eqglas}\pi_\nu(e^S)=(-1)^{b_1}\sum_{k_2,\dots,k_n}\tbinom{b_2+k_2}{b_2}\cdots\tbinom{b_N+k_N}{b_N}
				\pi_\nu(e^{S(k_2, \dots, k_N)}),
			\end{equation} where \[S(k_2,\dots,k_N)= \begin{matrix*}[l]
				1^{(a_1+a_2)}2^{(a_2-k_2)}  \cdots   N^{(a_N-k_N)} \\
				2^{(b_2+k_2)}  \cdots  N^{(b_N+k_N)}
			\end{matrix*}\] and the sum ranges over all  nonnegative integers $k_2,\dots,k_N$
			such that  $k_2+\dots+k_N=b_1 $ and $k_s \le a_s$ for all $s=2,\dots,N$.	\end{enumerate}	
	\end{lemma}
	We may think of the sum in the right hand side of eq. (\ref{eqglas}) as been taken over all ways of replacing the $b_1$ 1's in the second row of the tableaux $S$ with $k_2$ 2's, $k_3$ 3's, $\dots$, $k_N$ N's from the first row of $S$, where $k_2+ k_3+ \dots + k_N=b_1$. 
	
	Even though our paper \cite{MS3} concerns modular representations, the proof of the above lemma given there is valid for any field in place of $\mathbb{K}$ (in fact for any commutative ring). In \cite[Lemma 4.2]{MS3} we used the notation $\Delta_\mu$ for the Weyl module $K_\nu$. 
	\subsection{The duality functor $\Omega$}\label{omega}
Let us recall that there is an algebra involution on the ring of symmetric functions that sends the Schur function $s_\mu$ to $s_{\mu'}$ for every partition $\mu$ \cite[6.2]{F}. In terms of representations, we recall from \cite[p. 189]{AB} that, since the characteristic of $\mathbb{K}$ is zero, there is an involutory natural equivalence $\Omega$ from the category of polynomial representations of $G$ of degree $m$, where $N \ge m$, to itself that has the following properties.  \begin{enumerate} \item $\Omega(D(\alpha)) = \Lambda (\alpha)$ for all $\alpha \in \Lambda(N,m)$ and  $\Omega (K_\mu)=L_\mu$ for every $\mu \in \Lambda^+(N,r)$. More generally, if $\mu(1) \in \Lambda^{+}(N,r_1), \dots, \mu(q) \in \Lambda^{+}(N,r_q)$ are partitions such that $m_1 + \cdots + m_q =m$, then \[\Omega(K_{\mu(1)} \otimes \cdots \otimes K_{\mu(q)} ) = L_{\mu(1)} \otimes \cdots \otimes L_{\mu(q)}.\]\item The functor $\Omega$ preserves the comultiplication and multiplication maps of the Hopf algebras $D$ and $\Lambda$. 
	
	To be precise, this means that for all $(\alpha_1, \dots, \alpha_N) \in \Lambda(N,m)$ and all $s$ the images under $\Omega$ of the maps \begin{align*}&1\otimes \cdots \otimes \Delta_D \otimes \cdots \otimes 1 :D(\alpha_1, \dots, \alpha_s, \dots, \alpha_N) \to D(\alpha_1, \dots, \alpha'_s, \alpha''_s \dots, \alpha_N), \\&
		1\otimes \cdots \otimes \eta_D \otimes \cdots \otimes 1 :D(\alpha_1, \dots, \alpha_s, \alpha_{s+1} \dots, \alpha_N) \to D(\alpha_1, \dots, \alpha_s +\alpha_{s+1} \dots, \alpha_N)
	\end{align*}are the maps
	\begin{align*}&1\otimes \cdots \otimes \Delta_{\Lambda} \otimes \cdots \otimes 1 :\Lambda(\alpha_1, \dots, \alpha_s, \dots, \alpha_N) \to \Lambda(\alpha_1, \dots, \alpha'_s, \alpha''_s \dots, \alpha_N), \\&
		1\otimes \cdots \otimes \eta_{\Lambda} \otimes \cdots \otimes 1 :\Lambda(\alpha_1, \dots, \alpha_s, \alpha_{s+1} \dots, \alpha_N) \to \Lambda(\alpha_1, \dots, \alpha_s +\alpha_{s+1} \dots, \alpha_N)
	\end{align*}
	respectively. Here, $\alpha_s = \alpha'_s +\alpha''_s$ and $\Delta_D :D_{\alpha_s} \to D_{\alpha'_s} \otimes D_{\alpha''_s}$ 	and $\eta_D: D_{\alpha_s} \otimes D_{\alpha_{s+1}} \to D_{\alpha_s +\alpha_{s+1}}$
	are the indicated components of the comultiplication and multiplication maps of the divided power algebra $D$ respectively. Likewise, $\Delta_\Lambda :\Lambda^{\alpha_s} \to \Lambda^{\alpha'_s} \otimes \Lambda^{\alpha''_s}$ 	and $\eta_\Lambda: \Lambda^{\alpha_s} \otimes \Lambda^{\alpha_{s+1}} \to \Lambda^{\alpha_s +\alpha_{s+1}}$
	are the indicated components of the comultiplication and multiplication maps of the exterior algebra $\Lambda$ respectively.	
	
	\item We have $\Omega(\tau_{s,D}) = (-1)^{\alpha_s \alpha_{s+1}}\tau_{s,\Lambda}$ for all $\alpha =(\alpha_1, \dots, \alpha_N) \in \Lambda(N,m)$ and all $s$, where $\tau_{s,D}: D(\alpha) \to D(\alpha)$ (respectively, $\tau_{s,\Lambda}: \Lambda(\alpha) \to \Lambda(\alpha)$) is the map that interchanges the factors $D_{\alpha_s}$ and $D_{\alpha_{s+1}}$ (respectively, $\Lambda^{\alpha_s}$ and $\Lambda^{\alpha_{s+1}}$) and is the identity on the rest.\item The functor $\Omega$ is exact.\end{enumerate}
\subsection{Specht modules and the Schur functor $f$}	Suppose $N \ge m$. Recall from \cite{Gr} that the Schur functor is a functor $f$ from the category of homogeneous polynomial representations of  $G$ of degree $m$ to the category of $\mathfrak{S}_m$-modules. For $P$ an object in the first category,  $f(P)$ is the weight subspace $P_\alpha$ of $P$, where $\alpha = (1^m, 0^{N-m})$ and for $\theta :P \to Q$ a morphism in the first category, $f(\theta)$  is the restriction $P_\alpha \to Q_\alpha$ of $\theta$. 

We denote the space of column tabloids corresponding to a partition $\nu$ of $m$ by ${\tilde{M}}^\nu$, see \cite[Chapter 7.4]{F}. 

If $\lambda$ and $\mu$ are partitions such that $|\lambda|-|\mu| = m$, let us define the skew Specht $S^{\lambda / \mu }$ module by $S^{\lambda / \mu }:=f(L_{\lambda' / \mu'})$.

It is well known that $f$ is an exact functor. Me have  \begin{center}
	$f(L_{\nu})=S^{\nu'}$ and $f(\Lambda^{\nu}) = \tilde{M}^{\nu'}$.
\end{center}	\section{Relations of $\text{Lie}_n(m)$}
	The purpose of this section is to prove an identity in the multilinear component  $\mathrm{Lie}_n(m)$ of the free LAnKe that is crucial for our purposes. We know that the relations of the skew Weyl modules are parametrized by pairs of consecutive rows of the diagram of the skew partition, see Theorem \ref{skewmany}.   In Lemma \ref{relation} below, we prove an identity in $\mathrm{Lie}_n(m)$ for any pair of consecutive brackets.
	
	Let us denote $\mathrm{Lie}_n(m)$ by $\mathrm{Lie}_{n,k}$, where $m=(n-1)k+1$, when we want to emphasize the number of brackets $k$.
	
	The next lemma describes a particular relation in $\mathrm{Lie}_{n,3}$.
	
	We will denote a sequence of the form $a_1, \dots, a_{i-1}, a_{i+1}, \dots, a_q$ by $a_1, \dots, \widehat{a_{i}}, \dots, a_q$.
	
	\begin{lemma}\label{relationk=3} Let $k=3$ and $n \ge 2$. We have
		\begin{align*}
			(n-2)&[[z_0,z_1, \dots, z_{n-1}],w_1, \dots, w_{n-1}] \ + \\&[[z_0,w_1, \dots, w_{n-1}],z_1, \dots, z_{n-1}] \ + \\
			\sum_{i=1}^{n-1}(-1)^i\sum_{j=1}^{n-1}(-1)^{j+1}&[[z_0,w_j, z_1, \dots, \widehat{z_i}, \dots, z_{n-1}], z_i, w_1, \dots,  \widehat{w_j}, \dots, w_{n-1}] \\ = & \ 0
		\end{align*} for all $x_1, \dots, x_n, z_1, \dots, z_{n-1}, w_1, \dots, w_{n-1} \in [m]$, where $z_0=[x_1,\dots, x_n]$.
	\end{lemma}
	\begin{proof} From eq. (\ref{GJI}) we have \begin{align*}
			&[[z_0,z_1, \dots, z_{n-1}],w_1, \dots, w_{n-1}] \\\nonumber&= [[z_0,w_1, \dots, w_{n-1}],z_1, \dots, z_{n-1}] \ + \\\nonumber& \ \ \ \ \sum_{i=1}^{n-1}[z_0,z_1, \dots, z_{i-1}, [z_i,w_1, \dots, w_{n-1}], z_{i+1}, \dots, z_{n-1}].
		\end{align*}
		From the above equation and skew commutativity of the bracket, we obtain \begin{align}\label{relk311}
			&[[z_0,z_1, \dots, z_{n-1}],w_1, \dots, w_{n-1}]\\\nonumber&=[[z_0,w_1, \dots, w_{n-1}],z_1, \dots, z_{n-1}] \ + \\\nonumber&
			\ \ \ \ \sum_{i=1}^{n-1}(-1)^i[[z_i,w_1, \dots, w_{n-1}],z_0,z_1, \dots, \widehat{z_i}, \dots, z_{n-1}].
		\end{align}

		Let $A:= \sum_{i=1}^{n-1}(-1)^i[[z_i,w_1, \dots, w_{n-1}],z_0,z_1, \dots, \widehat{z_i}, \dots, z_{n-1}]$ be the last sum in (\ref{relk311}). To each summand \[[[z_i,w_1, \dots, w_{n-1}],z_0,z_1, \dots, \widehat{z_i}, \dots, z_{n-1}]\] of $A$ we apply eq. (\ref{GJI}) to obtain 
		\begin{align}\label{k31}A= &\sum_{i=1}^{n-1}(-1)^i[[z_i,z_0,z_1, \dots, \widehat{z_i}, \dots, z_{n-1}], w_1, \dots, w_{n-1}] \ + \\\nonumber &\sum_{i=1}^{n-1}(-1)^i\sum_{j=1}^{n-1}[z_i, w_{1}, \dots, w_{j-1},[w_j, z_0, z_1, \dots, \widehat{z_i}, \dots, z_{n-1}], w_{j+1}, \dots, w_{n-1}].
		\end{align}
		From skew commutativity of the bracket we have
		\begin{align*}\label{rel34}
			&[[z_i,z_0,z_1, \dots, \widehat{z_i}, \dots, z_{n-1}], w_1, \dots, w_{n-1}]\\&=(-1)^i[[z_0,z_1, \dots, z_{n-1}], w_1, \dots, w_{n-1}], \end{align*}
		and 
		\begin{align*}
&[z_i, w_{1}, \dots, w_{j-1},[w_j, z_0, z_1, \dots, \widehat{z_i}, \dots, z_{n-1}], w_{j+1}, \dots, w_{n-1}]\\&=(-1)^{j+1}[[z_0,w_j,z_1, \dots, \widehat{z_i},\dots, z_{n-1}], z_i, w_{1}, \dots, \widehat{w_j}\dots w_{n-1}]
		\end{align*}
		By substituting the last two equations in (\ref{k31}) we find
			\begin{align*}A= &(n-1)[[z_0,z_1, \dots, z_{n-1}],w_1, \dots, w_{n-1}] \ +\\\nonumber &\sum_{i=1}^{n-1}(-1)^i\sum_{j=1}^{n-1}(-1)^{j+1}[[z_0,w_j,z_1, \dots, \widehat{z_i},\dots, z_{n-1}], z_i, w_{1}, \dots, \widehat{w_j}\dots w_{n-1}].
		\end{align*}
		Substituting the last equation in eq. (\ref{relk311}) the desired result follows. \end{proof}	
	\begin{remark}\label{remk=3}The identity  of the previous lemma can be understood combinatorially. We observe that in the second term of the identity, all of the $z_1, \dots, z_{n-1}$ have been exchanged with the $w_1, \dots, w_{n-1}$. Also, in each term of the double sum, exactly one $z_i$ has been exchanged with one $w_j$. If $n=2$ the identity reads $0[[z_0,z_1],w_1]+[[z_0,w_1],z_1]-[[z_0,w_1],z_1]=0$.\end{remark} \begin{definition}A \textit{comb} of $s$ brackets on a set $X$ is a bracketed word of the form \begin{equation}\label{comb} [[\cdots[[x_0, x_{1,1} ,\dots, x_{1,n-1}], x_{2,1}, \dots, x_{2,n-1}], \dots, ], x_{s,1}, \dots, x_{s,n-1}], \end{equation} that has $s$ brackets,
	where $x_0, x_{i,j} \in X$.\end{definition} Combs are important for our purposes because of the following result of Friedmann, Hilton and Wachs \cite{FHW2}.
	\begin{proposition}[{\cite[Proposition 3.1]{FHW2}}]\label{propcombs}
		$\mathrm{Lie}_{n,k}$ is spanned by the combs of $k$ brackets on $[m]$.
	\end{proposition}
	
	The next lemma describes a specific relation between combs that depends on a choice of a pair of consecutive brackets. We number the brackets in a comb starting from the inner bracket to the outer. When we want to refer to a specific bracket we may denote this by putting a label on the bracket. For example when we label the ${s-1}$ bracket in the comb (\ref{comb}) we will write 
	\[ [[_{_{s-1}}\cdots[[x_0, x_{1,1} ,\dots, x_{1,n-1}], x_{2,1}, \dots, x_{2,n-1}], \dots, ], x_{s,1}, \dots, x_{s,n-1}].\]
	
	The main result of the present section is the following.
	\begin{lemma}\label{relation} Let $k \ge 3$ and $n \ge 2.$ Let $s$ be an integer such that $2 \le s \le k-1$. We have the following relations among combs in $\mathrm{Lie}_{n,k}$:
		\begin{align*}
			(n-2)&[[\cdots[[_{_s}z_0,z_1, \dots, z_{n-1}],w_1, \dots, w_{n-1}],x_{s+2, 1}, \dots,], x_{k,1}, \dots, x_{k,n-1}] \ + \\&[[\cdots[[_{_s}z_0,w_1, \dots, w_{n-1}],z_1, \dots, z_{n-1}],x_{s+2, 1},\dots,], x_{k,1}, \dots, x_{k,n-1}] \ + \\&
			\sum_{i=1}^{n-1}(-1)^i\sum_{j=1}^{n-1}(-1)^{j+1}[[\cdots[[_{_s}z_0,w_j, z_1, \dots, \widehat{z_i}, \dots, z_{n-1}],\\& z_i, w_1, \dots,  \widehat{w_j}, \dots, w_{n-1}],x_{s+2, 1}, \dots,], x_{k,1}, \dots, x_{k,n-1}] \\ = & \ 0
		\end{align*} for all combs $z_0$ with $s-1$ brackets on $[m]$ and for all $z_1, \dots, z_{n-1},$ $ w_1, \dots, w_{n-1}$, $x_{s+2, 1}, \dots, x_{s+2, n-1},$ $ \dots, x_{k, n-1} \in [m]$. 
	\end{lemma}
	\begin{proof}In the above formula, it is to be understood that all elements outside the $s+1$ bracket remain unchanged. Then the formula follows from Lemma \ref{relationk=3} and $n$-linearity of the bracket.
	\end{proof}	
	\begin{remark}\label{remk}Lemma \ref{relation} describes relations in $\mathrm{Lie}_{n,k}$ for any pair of consecutive brackets except for the pair (1,2). For each pair $(s,s+1)$ of consecutive brackets,  where $2 \le s \le k-1$, we have the situation of Remark \ref{remk=3}: In the second term of the equation of Lemma \ref{relation}, all of the $z_1, \dots, z_{n-1}$ of the $s$ bracket have been exchanged with the $w_1, \dots, w_{n-1}$ of the $s+1$ bracket. Also, in each term of the double sum, exactly one $z_i$ of the $s$ bracket has been exchanged with one $w_j$ of the $s+1$ bracket.\end{remark}	
	\section{The maps $\Phi_1$ and $\Phi_2$}
	\subsection{Strategy for the proof of the main result} In Sections 4-6 we prove that the answer to Question \ref{qu} is affirmative. In this subsection we discuss informally the main ideas of the proof. These will be illustrated for the case $k=4$.
	
	Recall we have the general linear group $G=GL_{N}(\mathbb{K})$ and the symmetric group $\mathfrak{S}_m$, $m=4n-3$. We assume $N \ge m$. Let $\lambda(4)$ be the partition $(n,n-1,n-1,n-1)$ and consider the space  of column tabloids $\tilde{M}^{\lambda(4)'}$ corresponding to the conjugate partition $\lambda(4)'$. As an $\mathfrak{S}_m$ module, $\tilde{M}^{\lambda(4)'}$ is isomorphic to the weight subspace of the tensor product $\Lambda^{n} \otimes \Lambda^{n-1} \otimes \Lambda^{n-1} \otimes \Lambda^{n-1}$ of exterior powers of the natural $G$-module of column vectors corresponding to the weight $(1^m, 0^{N-m})$. By skew commutativity of the bracket, we have an $\mathfrak{S}_m$-map 
	\[h:\tilde{M}^{\lambda(4)'} \to \mathrm{Lie}_{n,4}\] 
	that sends $x_1\cdots x_n \otimes y_1 \cdots y_{n-1}\otimes z_1 \cdots z_{n-1} \otimes w_1 \cdots w_{n-1} $ to the comb \[ [[[[x_1, \dots, x_n], y_1, \dots, y_{n-1}], z_1, \dots, z_{n-1}], w_1, \dots, w_{n-1}],\]
	for any permutation $x_1,\dots, x_n, y_1, \dots, y_{n-1}, z_1, \dots, z_{n-1}, w_1, \dots, w_{n-1}$ of $1,2, \dots, m$. By Proposition \ref{propcombs}, the map $h$ is surjective. 
	
	Now we define the skew partition $\alpha(n,4)=(n+2,n+1,n,n-1)/(2,2,1,0)$. The skew Specht module $S^{\alpha(n,4)'}$ corresponding to the conjugate skew partition is quotient of $\tilde{M}^{\lambda(4)'}$. Suppose for a moment that the relations of $S^{\alpha(n,4)'}$ are mapped to zero under the map $h$. Then we obtain a surjective map of $\mathfrak{S}_m$-modules $S^{\alpha(n,4)'} \to \mathrm{Lie}_{n,4}$. From the Littlewood-Richardson rule for skew Specht modules, it follows easily that every Specht module in the irreducible decomposition of $S^{\alpha(n,4)'}$ corresponds to a partition with at least 4 columns, if $n \ge 4$. Since the map $h$ is surjective, the same is true for $\mathrm{Lie}_{n,4}$. Hence with the notation of Question \ref{qu} we have $\gamma_{n,4} = 0$ as desired.
	
	Thus we must show that the relations of $S^{\alpha(n,4)'}$ are mapped to zero under $h$. In order to accomplish this, we find a new suitable description of $S^{\alpha(n,4)'}$ with the same generators but new relations. The fact that these new relations map to zero under $h:\tilde{M}^{\lambda(4)'} \to \mathrm{Lie}_{n,4}$ will be a consequence of Lemma \ref{relation}.
	
	In order to find a new description of $S^{\alpha(n,4)'}$, we will find a new description of the skew Weyl module $K_{\alpha(n,4)}$ with generators and relations, see Theorem \ref{mainweyl}(1). This is done by defining and studying suitable linear operators denoted $\Phi_1$ and $\Phi_2$ below.
	\subsection{The maps $\Phi_1, \Phi_2$}\label{phii}
	
	Let us begin by defining two maps that are crucial for Theorem \ref{mainweyl} below. \begin{definition}\label{mapsphi}Let $n \ge 2$. \begin{enumerate} \item Define $\Phi_1 \in \Hom _G(D(n,n-1), D(n,n-1))$ by \begin{align*}\Phi_1 &:= \phi_{S(1)}+ (-1)^n\phi_{S(2)}, \end{align*} where the tableaux $S(i) \in \mathrm{RSST}(n,n-1) $ are the following
			\[S(1):=\begin{matrix*}[l]
				1^{(n)}  \\
				2^{(n-1)}  \end{matrix*}, \ S(2):=\begin{matrix*}[l]
				12^{(n-1)}  \\
				1^{(n-1)} \end{matrix*}. \]
			\item Define $\Phi_2  \in \Hom _G(D(n-1,n-1), D(n-1,n-1))$ by \[\Phi_2 := (n-2)\phi_{S(3)}+ (-1)^{n-1}\phi_{S(4)} + \phi_{S(5)},\] where the tableaux $S(j) \in \mathrm{RSST}(n-1,n-1) $ are the following
			\[S(3):=\begin{matrix*}[l]
				
				1^{(n-1)} \\
				2^{(n-1)} \end{matrix*}, \ S(4):=\begin{matrix*}[l]
				
				2^{(n-1)} \\
				1^{(n-1)} \end{matrix*}, \ S(5):=\begin{matrix*}[l]
				
				1^{(n-2)}2 \\
				12^{(n-2)}\end{matrix*}. \]
		\end{enumerate}
	\end{definition}
	\begin{remark} A map closely related  to $\Phi_1$ was defined and studied in \cite[Section 3.1]{MS7} (the notation $\gamma_1$ was used there). The motivation for the definition of $\Phi_2$ will be clear in the proof of Theorem \ref{main} (cf. eq. (\ref{omegaphi3})). Roughly speaking, it turns out that when we apply to the map $\Phi_2$ the duality functor $\Omega$ followed by the Schur functor, the image of the resulting map is given by the relations of Lemma \ref{relation}.
	\end{remark}
	\begin{remark}\label{cyc} It well known that for any sequence $\alpha=(\alpha_1, \dots, \alpha_r)$ of nonnegative integers, where $N \ge r$, the $G$-module $D(\alpha)$ is cyclic and a generator is the element $1^{(\alpha_1)}\otimes \cdots \otimes r^{(\alpha_r)}$.
\end{remark}
	
	We want to describe the  image of the map $\Phi_2$.
	Recall from Definition \ref{thetat} that we have the map \[\theta_2 : D(n+1,n-3) \to D(n-1,n-1)\] that satisfies  \begin{equation}\label{coktheta2}\coker{(\theta_2)} \simeq K_\xi,\end{equation} where $\xi$ is the skew shape $\xi := (n,n-1)/(1).$
	\begin{lemma}\label{newim}Let $n \ge 3$. Then, with the above notation, the images of the maps 
		\begin{align*}
			&\Phi_2: D(n-1,n-1) \to D(n-1,n-1),\\&
			\theta_2 : D(n+1,n-1) \to D(n-1,n-1)
		\end{align*}
	are equal.
	\end{lemma}
	\begin{proof}
		First we will determine a new description of $\Ima(\Phi_2)$. This will be done using ideas from \cite[Section 3.1]{MMS}, but our presentation is self contained.
		
		Recall that we have the decomposition into irreducible $G$-modules \begin{equation}\label{pierin-1}D(n-1,n-1) = \bigoplus_{0 \le t \le n-1}K_{(n-1+t,n-1-t)}\end{equation} by Pieri's rule \cite[(2.3.5) Corollary]{W}. In order to have concise notation, let $\nu (t)$ be the partition \[\nu (t)=(n-1+t, n-1-t).\]We have the projection \begin{equation}\label{pi_T}\pi_T:D(n-1, n-1) \xrightarrow{\phi_T} D(\nu (t)) \xrightarrow{\pi_{\nu (t)}} K_{\nu (t)}\end{equation} onto the irreducible $K_{\nu(t)}$, where $T$ is the tableau \[T:=\begin{matrix*}[l]
			1^{(n-1)} 2^{(t)} \\
			2^{(n-1-t)} \end{matrix*}\] (cf. Definition \ref{phis}). 
			
			We will find the decomposition into irreducible $G$-modules of the image of the map $\Phi_2$ by examining the compositions  \begin{equation}\label{pi_TPhi_2}D(n-1,n-1) \xrightarrow{\Phi_2} D(n-1,n-1) \xrightarrow{\pi_T} K_{\nu (t)}\end{equation} for $t=0, \dots, n-1$.  This will be done by  determining the $t \in \{ 0,1, \dots, n-1 \}$ such that the composition (\ref{pi_TPhi_2}) is nonzero.
			
			\textit{Claim}.  The composition (\ref{pi_TPhi_2}) is nonzero if and only if $2 \le t \le n-1$.
		
		We compute the image of the element $1^{(n-1)} \otimes 2^{(n-1)}$ under the map (\ref{pi_TPhi_2}).  From Definition \ref{mapsphi} we have  \begin{align}\label{Phi3*}\nonumber\Phi_2(1^{(n-1)} \otimes 2^{(n-1)}) =&(n-2) 1^{(n-1)}\otimes 2^{(n-1)}+(-1)^{n-1} 2^{(n-1)} \otimes 1^{(n-1)}  \\&+1^{(n-2)}2\otimes 12^{n-2}.\end{align}
		Using (\ref{Phi3*}) we obtain
		\begin{align}\label{piPhi2}\nonumber\pi_T \circ \Phi_2(1^{(n-1)} \otimes 2^{(n-1)}) =&(n-2) \pi_T(1^{(n-1)}\otimes 2^{(n-1)})\\\nonumber&+(-1)^{n-1} \pi_T(2^{(n-1)} \otimes 1^{(n-1)})  \\&+\pi_T(1^{(n-2)}2\otimes 12^{n-2}).\end{align}
		Using the definition of the map $\pi_T$ given in  (\ref{pi_T}) we compute 
		\begin{align}
			&\label{4.6}\pi_T(1^{(n-1)}\otimes 2^{(n-1)})= \pi_{\nu (t)}(1^{(n-1)}2^{(t)}\otimes 2^{(n-1-t)}), \\\label{4.7}
			&\pi_T(2^{(n-1)} \otimes 1^{(n-1)})= \pi_{\nu (t)}(1^{(t)}2^{(n-1)}\otimes 1^{(n-1-t)}), \\\label{4.8}
			&\pi_T(1^{(n-2)}2\otimes 12^{(n-2)}) = (n-1)t\pi_{\nu (t)}(1^{(n-1)}2^{(t)}\otimes 2^{(n-1-t)}) + \\\nonumber  & \ \ \ \ \ \ \ \ \ \ \ \ \ \ \ \ \ \ \ \ \ \ \ \ \ \ \  \ \ \ \  \ \ \ \  (t+1)\pi_{\nu (t)}(1^{(n-2)}2^{(t+1)}\otimes 12^{(n-2-t)}).
		\end{align}
		Now we apply Lemma \ref{lemglas}(2) to the right hand side of (\ref{4.7}) and to the second summand in the right hand side of (\ref{4.8}) obtaining respectively
		\begin{align}\label{4.9}&\pi_{\nu (t)}(1^{(t)}2^{(n-1)} \otimes 1^{(n-1-t)})= (-1)^{n-1-t}\pi_{\nu (t)}(1^{(n-1)} 2^{(t)}\otimes 2^{(n-1-t)}), \\\label{4.10} &\pi_{\nu (t)}(1^{(n-2)}2^{(t+1)}\otimes 12^{n-2-t})=-\tbinom{n-1-t}{1}\pi_{\nu (t)}(1^{(n-1)}2^{(t)}\otimes 2^{(n-1-t)}).\end{align}
		Substituting equations (\ref{4.6}) - (\ref{4.10}) in (\ref{piPhi2}) we obtain
		\begin{equation}\label{piPhi3*2}\pi_T \circ \Phi_2(1^{(n-1)} \otimes 2^{(n-1)}) =q\pi_{\nu (t)}(1^{(n-1)}2^{(t)}\otimes 2^{(n-1-t)}),\end{equation}
		where $q:=n-2 + (-1)^t +(n-1)t-(t+1)(n-1-t)$. Simplifying we have \[q=(-1)^t+t^2+t-1\]
		and we see  that \begin{equation}\label{q}q=0 \Leftrightarrow t\in \{0,1\}.\end{equation} 
		
		Let us consider eq. (\ref{piPhi3*2}). The element $\pi_{\nu (t)}(1^{(n-1)}2^{(t)}\otimes 2^{(n-1-t)})$ in the right hand side of this equation is a semistandard basis element of the Weyl module $K_{\nu(t)}$ and therefore is nonzero. On the other hand, the element $1^{(n-1)} \otimes 2^{(n-1)}$ is a generator of the $G$-module $D(n-1,n-1)$ according to Remark \ref{cyc}. Thus from eq. (\ref{piPhi3*2}) and eq. (\ref{q}) we conclude that the composition (\ref{pi_TPhi_2}) is zero if and only if $t \in \{0,1\}$. In other words, the composition (\ref{pi_TPhi_2}) is nonzero if and only if $2 \le t \le n-1$. This proves the claim.
		
		Now since the irreducible decomposition (\ref{pierin-1}) is multiplicity free, the Claim implies  that \[\Ima(\Phi_2) = \bigoplus_{2 \le t \le n-1}K_{(n-1+t, n-1-t)}.\] In other words, the cokernel of the map $\Phi_2$ is  \[\coker(\Phi_2)= K_{(n-1,n-1)} \oplus K_{(n,n-2)}.\]
		However, $K_{(n-1,n-1)} \oplus K_{(n,n-2)} = K_\xi$ where $\xi$ is the skew partition $\xi= (n,n-1) / (1),$
		according to Pieri's rule for skew shapes, see \cite[(2.3.7) Corollary]{W}. 
		
		On the other hand, we know from (\ref{coktheta2}) that $\coker(\theta_2)=K_\xi$. Now since the irreducible decomposition (\ref{pierin-1}) is multiplicity free, we conclude that $\Ima(\Phi_2) = \Ima(\theta_2)$.\end{proof} 
	
	\section{Main result for Weyl modules}
	In this section we prove our main result concerning certain skew Weyl modules, Theorem \ref{mainweyl}. We begin with the relevant definitions.
	\subsection{Definitions and statement of main result}
	\begin{definition}
		Let $k\ge 2$ be an integer. We define the partition \[\lambda(k):=(n, (n-1)^{k-1}),\] where $(n-1)^{k-1}$ means that the part $n-1$ appears $k-1$ times.
	\end{definition}
	
	Recall we have the maps $ \Phi_1 : D(n,n-1) \to D(n,n-1)$ and $ \Phi_2 : D(n-1,n-1) \to D(n-1,n-1)$ of Definition \ref{mapsphi}. Using these we give the following definition.
	\begin{definition}\label{defcok} For each $s=1, \dots, k-1$, let $M:=M_s=D(\lambda(k))$. Define the $G$-module $U(n,k)$ as the cokernel of the map 
		\[  \bigoplus_{s=1}^{k-1} M_s \xrightarrow{\Phi(n,k)} M,\]
		where \begin{itemize}
			\item \ the restriction of the map $\Phi(n,k)$ to $M_1$ is the map \[D(n,n-1) \otimes D({(n-1)^{k-2}}) \xrightarrow{\Phi_1\otimes 1} D(n,n-1) \otimes D({(n-1)^{k-2}}),\] where $1$ is the identity map on $D({(n-1)^{k-2}}),$ and
			\item \ the restriction of the map $\Phi(n,k)$ to $M_s$, $s>1$, is the map \begin{align*}&D(n,(n-1)^{s-1}) \otimes D(n-1, n-1) \otimes D((n-1)^{k-1-s}) \xrightarrow{1 \otimes \Phi_2\otimes 1} \\&D(n,(n-1)^{k-1-s}) \otimes D(n-1, n-1) \otimes D(n,(n-1)^{s-1}),\end{align*} where $1$ is the appropriate identity map.
	\end{itemize} \end{definition}
	
	\begin{example} Let $k=4$. Then \[M=D(n,n-1,n-1,n-1)=D_n \otimes D_{n-1} \otimes D_{n-1} \otimes D_{n-1}\] and $U(n,k)$ is by definition the cokernel of the map \[M \oplus M \oplus M \xrightarrow{\Phi_1 \otimes 1 \otimes 1 \ + \ 1 \otimes \Phi_2 \otimes 1 \ + \ 1 \otimes 1\otimes \Phi_2 } M.\]
	\end{example}
	\begin{definition}\label{alphank} For integers $n,k \ge 2$, let $\alpha(n,k)$ be the skew partition \[\alpha(n,k) := (n+k-2, n+k-3, \dots, n,n-1)/(k-2, k-2, k-3, k-4, \dots, 2,1,0).\]
	\end{definition}	
	For example, if $n=k=5$ then the Young diagram of $\alpha(n,k)$ looks like 
	
	\begin{center}
		\ytableausetup{boxsize=1.2em}
		\ydiagram{3+5,3+4,2+4,1+4,4}
	\end{center}
	
	We remark that $\alpha(n,k)$ is the conjugate skew partition of $\beta{(n,k)}$ defined in (\ref{beta}).
	
	We note that the first row of the Young diagram of 	$\alpha(n,k)$ has $n$ cells and every other row has $n-1$ cells. Also, the first two rows have $n-1$ overlaps and every other pair of consecutive rows has $n-2$ overlaps. In particular, $\alpha(n,k) = \lambda(k)$. Hence from Theorem \ref{skewmany}  we get the following corollary.
	\begin{corollary}\label{coro}As $G$-modules, the skew Weyl module $K_{\alpha(n,k)}$ is isomorphic to the cokernel of the map \begin{align}&\theta_{\alpha(n,k)} : A_1 \oplus A_2 \oplus \dots \oplus A_{k-1}\to D(\lambda(k)),\\\nonumber&\theta_{\alpha(n,k)} = \theta_1 \otimes 1  + 1 \otimes \theta_2 \otimes 1   + \dots + 1  \otimes \theta_2,\end{align}
		where \begin{align*}&A_1:=D(n+1,n-2) \otimes D((n-1)^{k-2}),
			\\&A_2:=D(n)\otimes D(n+1, n-3) \otimes D((n-1)^{k-3}),
			\\& A_3 :=D(n,n-1)\otimes D(n+1, n-3) \otimes D((n-1)^{k-4}),
			\\&\ \ \ \ \ \vdots 
			\\&A_{k-1}:=D(n,{(n-1)}^{k-3})\otimes D(n+1,n-3),
		\end{align*}
		and the maps are \begin{align*}&\theta_1 \otimes 1 : A_1 \to D(\lambda (k)), \\ &1 \otimes \theta_2 \otimes 1 : A_2 \to D(\lambda (k)), \\ &1 \otimes \theta_2 \otimes 1 : A_3 \to D(\lambda (k)), \\ & \ \ \ \ \ \vdots \\ &1\otimes \theta_2 : A_{k-1} \to D(\lambda (k)).\end{align*} \end{corollary}
	
	A simple yet important observation is the following.
	\begin{remark}\label{column}
		If $n \ge k$, then the Young diagram of the skew partition $\alpha(n,k)$ has a column of length $k$. 
	\end{remark}	
	
	Our main result for Weyl modules is the following.
	\begin{theorem}\label{mainweyl}
		Suppose $N \ge (n-1)k+1$. \begin{enumerate}
			\item[\textup{(1)}] The $G$-modules $U(n,k)$ and $K_{\alpha(n,k)}$ are isomorphic.
			\item[\textup{(2)}]If $n \ge k$, then for every partition $\nu$ of length at most $k-1$, the  multiplicity of the irreducible module $K_\nu$ in $U(n,k)$ is equal to zero.
		\end{enumerate}
	\end{theorem}

	\subsection{The cokernel of the map $\Phi_1: D(n,n-1) \to D(n,n-1)$}
		In preparation for the proof of Theorem \ref{mainweyl}, we will recall a result from \cite{MMS}.
	\begin{definition}\cite[Defintion 5.1]{MMS} \begin{enumerate}\item Let $\beta_{n-1}$ be the map of tensor product of exterior powers \[\beta_{n-1}: \Lambda^n \otimes \Lambda^{n-1} \to \Lambda^n \otimes \Lambda^{n-1}\] given by the composition
			\[  \Lambda^n \otimes \Lambda^{n-1} \xrightarrow {\Delta \otimes 1} \Lambda^1 \otimes \Lambda^{n-1} \otimes \Lambda^{n-1} \xrightarrow{\tau}  \Lambda^1 \otimes \Lambda^{n-1} \otimes \Lambda^{n-1} \xrightarrow{\eta \otimes 1}  \Lambda^{n} \otimes \Lambda^{n-1}, \] 
			where $\Delta: \Lambda^{n} \to \Lambda^1 \otimes \Lambda^{n-1}$ (respectively, $\eta: \Lambda^{1} \otimes \Lambda^{n-1} \to \Lambda^{n}$) is the indicated component of the comultiplication map (respectively, multiplication map) of the exterior algebra $\Lambda$, and $\tau: \Lambda^{n-1} \otimes  \Lambda^{n-1} \to \Lambda^{n-1} \otimes  \Lambda^{n-1}$ is the defined by $\tau(x \otimes y) =y \otimes x$, for $x, y \in \Lambda^{n-1} $. \item Let $\gamma_{n-1}$ be the map \[\gamma_{n-1}: \Lambda^n \otimes \Lambda^{n-1} \to \Lambda^n \otimes \Lambda^{n-1}, \ \gamma_{n-1}(x\otimes y) =x \otimes y - \beta_{n-1}(x \otimes y). \] \end{enumerate}\end{definition}
	To be precise, the maps $\beta_{n-1}$ and $\gamma_{n-1}$ above are the special cases of the maps $\beta_k$ and $\gamma_k$ defined in \cite[Definition 5.1]{MMS} for $a=n$ and $ b=k=n-1$ in the notation of loc. cit. The maps $\beta_{n-1}$ and $\gamma_{n-1}$ should not be confused with the symmetric group modules $\beta_{n,k}$ and $\gamma_{n,k}$ of the Introduction.
	
	We recall the following special case of \cite[Corollary 5.4]{MMS}. Here $L_{(n,n-1)}$ denotes the Schur module corresponding to the partition $(n, n-1)$, see Section 2.2.
	\begin{lemma}\label{l412} Suppose $N \ge 2n-1$. Then $\coker(\gamma_{n-1}) \simeq L_{(n,n-1)}$.
	\end{lemma}
	We note without pursuing details that the previous Lemma also follows from \cite[Theorem 1.3]{FHSW} by applying first the inverse Schur functor and then the functor $\Omega$.
	
	Now we may identify the cokernel of the map $\Phi_1: D(\lambda) \to D(\lambda)$ of Definition \ref{mapsphi}(1).
	\begin{corollary}\label{c413} Suppose $N \ge 2n-1$. Then $\coker(\Phi_1) \simeq K_{(n,n-1)}$.
	\end{corollary}
	\begin{proof}
		Consider the functor $\Omega$ of Section \ref{omega}. Using the properties (1) - (3) of $\Omega$, it follows that the image of the map  \[\gamma_{n-1}: \Lambda^n \otimes \Lambda^{n-1} \to \Lambda^n \otimes \Lambda^{n-1}\] under $\Omega$ is the map $\Phi_1: D(n,n-1) \to D(n,n-1)$. Since $\Omega$ is an exact functor, we have  \[\coker{(\Phi_1)} \simeq \Omega (\coker{\gamma_{n-1}}).\]
		From Lemma \ref{l412} we have \[\Omega (\coker{\gamma_{n-1}}) \simeq \Omega (L_{(n,n-1)})\]
		and according to property (1) of $\Omega$,\[\Omega (L_{(n,n-1)}) \simeq K_{(n,n-1)}.\] Hence $\coker(\Phi_1) \simeq K_{(n,n-1)}$.\end{proof}
	
\subsection{Proof of Theorem \ref{mainweyl}} We proceed with the proof of Theorem \ref{mainweyl}. Suppose $N \ge (n-1)k+1$. Recall the $G$-module $U(n,k)$ from Definition \ref{defcok}.
	\begin{proof} (1) Recall we have the maps \begin{align*}
			&\Phi_1 :D(n,n-1) \to D(n,n-1),\\
			&\theta_1 :D(n+1,n-2) \to D(n,n-1)
		\end{align*}
		given in Definition \ref{mapsphi}(1) and in (\ref{skew2}) respectively.
		From Corollary \ref{c413} we have \begin{equation}\label{cokPhi1}
			\coker(\Phi_1) \simeq K_{(n,n-1)}.
		\end{equation} From (\ref{skew2}) we have \begin{equation}\label{coktheta1}\coker{(\theta_1)} \simeq K_{(n,n-1)}.\end{equation} Since the decomposition of the codomain $D(n,n-1)$ of these maps into a sum of irreducible $G$-modules is multiplicity free according to Pieri's rule, we conclude from (\ref{cokPhi1}) and (\ref{coktheta1}) that \begin{equation}\label{theta1}
			\Ima(\Phi_1) = \Ima(\theta_1). \end{equation}
		According to Lemma \ref{newim}, 
		\begin{equation}\label{theta2}
			\Ima(\Phi_2) = \Ima(\theta_2). \end{equation}
		From (\ref{theta1}) and (\ref{theta2}) we conclude that the maps 	\[  \bigoplus_{s=1}^{k-1} M_s \xrightarrow{\Phi(n,k)} M,\] of Definition \ref{defcok} and \[\theta_{\alpha(n,k)} : A_1 \oplus A_2 \oplus \dots \oplus A_{k-1}\to D(\lambda(k))\] of Corollary \ref{coro} have equal images. Hence $U(n,k) \simeq K_{\alpha(n,k)}.$
		
		(2) The decomposition of $K_{\alpha(n,k)}$ into irreducibles is given by the Littlewood-Richardson rule \cite[Section 5.2]{F}: The multiplicity of $K_\nu$ in $K_{\alpha(n,k)}$ is equal to the number of Littlewood-Richardson semistandard tableaux of shape $\alpha(n,k)$ and weight $\nu$. Assume $n \ge k$. By Remark \ref{column} every  semistandard tableau $T$ of shape $\alpha(n,k)$ contains a column of length $k$ with entries strictly increasing from top to bottom. Hence, the weight of $T$ is a partition of length at least $k$. Thus, if $\nu$ has at most $k-1$ parts, no such $T$ can exist and the multiplicity of $K_\nu$ in $K_{\alpha(n,k)}$ is equal to zero. \end{proof}
	\section{Main result for LAnKes} Let us recall the setup and notation from the Introduction. Let $\mathrm{Lie}_n(m)$ be the multilinear component of the free LAnKe on the set $[m]$, where $m= (n-1)k+1$. We have the action of the symmetric group $\mathfrak{S}_m$ on $\mathrm{Lie}_n(m)$ described in the Introduction and we denote the corresponding representation by $\rho_{n,k}$. According to Theorem \ref{FHWdec}, we have  $\rho_{n,k}=\beta_{n,k} \oplus \gamma_{n,k}$ as $\mathfrak{S}_m$-modules, where the decomposition of  $\beta_{n,k}$ into irreducibles is obtained by adding a row of length $k$ to the top of each Young diagram in the decomposition of $\rho_{n-1,k}$, and $\gamma_{n,k}$ is a $\mathfrak{S}_m$-module all of whose irreducibles have Young diagrams with at most $k-1$ columns. 
	
	We have the skew partition $\alpha(n,k)$ of $m$ given in Definition \ref{alphank} and the skew Specht module $S^{\alpha(n,k)'}$  corresponding to the conjugate skew shape $\alpha(n,k)'$ of $\alpha(n,k)$.

	The main result of this paper is the following.
	
	\begin{theorem}\label{main} Suppose $N \ge m$. \begin{enumerate}
			\item[\textup{(1)}] There is a surjective map of  $\mathfrak{S}_m$-modules $S^{a(n,k)'} \to \mathrm{Lie}_n(m)$.
			\item[\textup{(2)}] If $n \ge k$, then $\gamma_{n,k} =0$.
		\end{enumerate}
	\end{theorem}
	\begin{proof}We want to use Theorem \ref{mainweyl} to obtain a result about Specht modules of the symmetric group $\mathfrak{S}_m$. In other words we want to apply the functor $\Omega$ and then the Schur functor $f$ (see Sections 2.5 and 2.6). 
		
		Recall we have the partition $\lambda(k)=(n, (n-1)^{k-1})$ of $m$. Its conjugate is $\lambda (k)'=(k^{n-1},1)$. Also we have the space of column tabloids $\tilde{M}^{\lambda(k)'}$ corresponding to the partition $\lambda(k)'$.
		
		(1) From skew commutativity of the bracket in the definition of LAnKe given in the Introduction, it follows that there is  map of $\mathfrak{S}_m$-modules 
		\begin{align*}h: \ &\tilde{M}^{\lambda(k)'} \to \text{Lie}_n(m), \\ &|x_{1,0} x_{1,1} \cdots x_{1,n}|x_{2,1} \cdots x_{2,n-1}| \cdots |x_{k,1} \cdots x_{k,n-1}| \mapsto \\& [[\cdots[[x_{1,0}, x_{1,1}, \dots, x_{1,n}],x_{2,1}, \dots, x_{2,n-1}], \dots,],x_{k,1}, \dots, x_{k,n-1}],  \end{align*}
		where $x_{i,j} \in [m]$ and  $|x_{1,0} x_{1,1} \cdots x_{1,n}|x_{2,1} \cdots x_{2,n-1}| \cdots |x_{k,1} \cdots x_{k,n-1}|$ denotes the column tabloid corresponding to the partition $\lambda (k)'=(k^{n-1},1)$ that has first column $x_{1,0} x_{1,1} \cdots x_{1,n}$, second column $x_{2,1} \cdots x_{2,n-1}$ etc. According to Lemma \ref{propcombs}, $\mathrm{Lie}_n(m)$ is generated by combs of $k$ brackets on $[m]$. Hence the map $h$ is surjective.
		
		Applying the functor $\Omega$ to the map $\Phi(n,k)$ of Definition \ref{defcok} yields the exact sequence of $G$-modules
		\begin{equation}\label{omegacok}\bigoplus_{s=1}^{k-1}\Omega(M_s) \xrightarrow{\Omega(\Phi(n,k))} \Omega(M) \to \Omega(U(n,k)) \to 0,
		\end{equation}
		where $\Omega(M_s)=\Omega(M)=\Lambda(\lambda(k))$ for $s=1, \dots, k-1$.
		It is straightforward to verify from the definition of the maps $\Phi_1, \Phi_2$ and properties (3) and (4) of the functor $\Omega$ (Section \ref{omega}), that the maps \begin{align*}
			&\Omega(\Phi_1):\Lambda(n,n-1) \to \Lambda(n,n-1)
			\\&\Omega(\Phi_2):\Lambda(n-1,n-1) \to \Lambda(n-1,n-1)\end{align*} are given as follows
		\begin{equation}\label{omegaphi1} \Omega(\Phi_1)(x \otimes y)= x \otimes y -\sum_{i=1}^{n}(-1)^{i-1}x_iy \otimes x_1 \cdots \widehat{x_i} \cdots x_{n},
		\end{equation}
		for all $x=x_1\cdots x_{n}$ and $y =y_1 \cdots y_{n-1}$, where  $x_i,y_j \in \Lambda^1$, and
		\begin{align}\label{omegaphi3}\Omega(\Phi_2)&(z \otimes w)= (n-2)z \otimes w + w \otimes z + \\\nonumber&
			\sum_{i=1}^{n-1}(-1)^i\sum_{j=1}^{n-1}(-1)^{j+1} w_j z_1 \cdots \widehat{z_i} \cdots z_{n-1}\otimes z_i w_1 \cdots \widehat{w_j} \cdots w_{n-1},
		\end{align}
		for all $z=z_1\cdots z_{n-1}$ and $w =w_1 \cdots w_{n-1}$ where  $z_i,w_j \in \Lambda^1$.
		
		Applying the Schur functor to the exact sequence (\ref{omegacok}) we obtain the exact sequence 
		\begin{equation}\label{fomegacoker}\bigoplus_{s=1}^{k-1}\tilde{M}^{\lambda(k)'}  \xrightarrow{f(\Omega({\Phi(n,k)})} \tilde{M}^{\lambda(k)'} \to f(\Omega(U(n,k))) \to 0.
		\end{equation}
		
		Consider $s=1$ in (\ref{fomegacoker}). From Definition \ref{defcok} it follows that the restriction of the map $f(\Omega({\Phi(n,k)})$ to the $s=1$ summand is
		\begin{equation}\label{rest1}\tilde{M}^{\lambda(k)'} \xrightarrow{f(\Omega({\Phi_1}))\otimes 1} \tilde{M}^{\lambda(k)'}. \end{equation}
		From (\ref{omegaphi1}) and the generalized Jacobi identity (\ref{GJI}), it follows that the image of the map (\ref{rest1}) is contained in the kernel of $h$. 
		
		Consider $s>1$ in (\ref{fomegacoker}). From Definition \ref{defcok} it follows that the restriction of the map $f(\Omega({\Phi(n,k)})$ to the $s>1$ summand is
		\begin{equation}\label{rest2}\tilde{M}^{\lambda(k)'} \xrightarrow{1 \otimes f(\Omega({\Phi_2}))\otimes 1} \tilde{M}^{\lambda(k)'}. \end{equation}From eq. (\ref{omegaphi3}) and Lemma \ref{relation} it follows that the image of the map (\ref{rest2}) is contained in the kernel of $h$. 
		
		Hence we conclude from (\ref{fomegacoker}) that $h$ induces a surjective map of $\mathfrak{S}_m$-modules
		\[f(\Omega(U(n,k))) \to \mathrm{Lie}_n(m).\]
		
		From Theorem \ref{mainweyl}(1) we have
		\[f(\Omega(U(n,k))) \simeq f(\Omega(K_{\alpha(n,k)}))\] as $\mathfrak{S}_m$-modules. From property (1) of the functor $\Omega$ we obtain \[f(\Omega(K_{\alpha(n,k)})) \simeq f(L_{\alpha(n,k)}) \] as $\mathfrak{S}_m$-modules. But for the Schur functor $f$ we have \[f(L_{\alpha(n,k)}) \simeq S^{\alpha(n,k)'} \] as $\mathfrak{S}_m$-modules. Thus we conclude that there is a surjective map of $\mathfrak{S}_m$-modules $S^{a(n,k)'} \to \mathrm{Lie}_n(m)$.
		
(2) We saw in the previous paragraph that $S^{\alpha(n,k)'} \simeq f(\Omega(U{(n,k)}))$. Suppose $n \ge k$. From Theorem \ref{mainweyl}(2) we know that the multiplicity in $U(n,k)$ of any $K_\nu$, where $\nu$ is a partition of $m$ of length at most $k-1$, is equal to zero. Since the functors $\Omega$ and $f$ satisfy \[\Omega(K_\nu) =L_\nu \ \ \text{and} \ \  f(L_\nu)=S^{\nu'}\] for all partitions $\nu$ of $m$, we conclude that the multiplicity in $f(\Omega(U{(n,k)}))$ of any Specht module whose diagram has at most $k-1$ columns, is equal to zero. Thus the same is true $S^{\alpha(n,k)'}$. From part (1) of the present theorem, the same is true for $\mathrm{Lie}_n(m)$. This means that $\gamma_{n,k}=0$.\end{proof}
	
	As mentioned in the Introduction, the decomposition into irreducibles of $\rho_{n,k}$ is known in the following cases
	\begin{itemize} \item \  $n=2$ (\cite{KW}),
		\item \ $k=2$ (\cite{FHSW}),
		\item  \ $k=3$ (\cite{FHW2, MS7}),
		\item \ $n=3$ and $k=4$ (\cite{FHW2}).
	\end{itemize}
	
	From the previous theorem we may deduce the decomposition into irreducibles of $\rho_{n,4}$ for all $n \ge 3$.\begin{corollary}\label{cordec} Let $n\ge3$. Then the decomposition into irreducibles of $\rho_{n,4}$ is
		\begin{align*}\rho_{n,4} = &S^{(4^{n-1}1)} \oplus S^{(4^{n-2}32)} \oplus S^{(4^{n-2}31^2)} \oplus S^{(4^{n-2}2^21)} \oplus \\&S^{(4^{n-2}21^3)} \oplus S^{(4^{n-3}3^21^3)} \oplus S^{(4^{n-3}32^3)}.\end{align*}
	\end{corollary}
	\begin{proof}
		From Theorem \ref{main} and Theorem \ref{dec} it follows that the decomposition of $\rho_{n,4}$ into a sum of irreducibles is obtained by adding $n-3$ rows of length $4$ to the Young diagram of each irreducible summand of $\rho_{3,4}$. The decomposition of $\rho_{3,4}$ into a sum of irreducibles was obtained in Theorem 4.3 of version 2 of \cite{FHW2}. From this the result follows.\end{proof}
	
	It was proved in \cite{FHW2} that for $n <k$ it is possible that $\gamma_{n,k}$ is nonzero. Very little seems to known about the irreducible decomposition of $\gamma_{n,k}$ in general. The first part of Theorem \ref{main} implies the following upper bound on the multiplicities of irreducibles in $\rho_{n,k}$ and hence in $\gamma_{n,k}$.
	\begin{corollary}\label{cordecsp} Let $n\ge3$ and $\nu$ a partition of $m$. Then, the multiplicity of $S^{\nu}$ in $\rho_{n,k}$ is less than or equal to the Littlewood-Richardson coefficient $c^{\lambda}_{\mu, \nu}$, where \begin{align*}\lambda &= (k^{n-1}, k-1, k-2, \dots, 1), \\ \mu&=(k-1, k-2, \dots, 2).
	\end{align*}\end{corollary}

	\section*{Acknowledgments} Dimitra-Dionysia Stergiopoulou was supported by a postdoctoral scholarship from the Christos Papakyriakopoulos Bequest of the National Technical University of Athens.

	\end{document}